\documentclass[11pt]{article} 
  
\usepackage{amsfonts}
\usepackage{amssymb,amsmath} 

\newcommand{\nicefrac}[2]
{\leavevmode \kern.1em\raise.5ex\hbox{\the\scriptfont0 #1}
             \kern-.1em/\kern-.15em\lower.25ex
             \hbox{\the\scriptfont0 #2}}

\newtheorem{Theo}{Theorem}{\alph{enumi}}
\newenvironment{theorem}{\begin{Theo}\hspace{-0.2cm}: }{\end{Theo}}

\newtheorem{Pro}{Proposition}{\alph{enumi}}

\newtheorem{Co}{Corollary}{\alph{enumi}}
\newenvironment{corollary}{\begin{Co}\hspace{-0.2cm}: }{\end{Co}}

\newtheorem{As}{Assumption}{\alph{enumi}}

\newtheorem{Le}{Lemma}{\alph{enumi}}
\newenvironment{lemma}{\begin{Le}\hspace{-0.2cm}: }{\end{Le}}

\newtheorem{Fo}{Folgerung}{\alph{enumi}}

\newtheorem{De}{Definition}{\alph{enumi}}
\newenvironment{definition}{\begin{De}\hspace{-0.2cm}:\rm }{\end{De}}

\newtheorem{Be}{Bemerkung}{\alph{enumi}}

\newtheorem{Ex}{Example}{\alph{enumi}}
\newenvironment{example}{\begin{Ex}\hspace{-0.2cm}:\rm }{\end{Ex}}

\newcommand{\ol}{\overline}
\newcommand{\bee}{\begin{equation}}
\newcommand{\ee}{\end{equation}}

\newcommand{\bea}{\begin{eqnarray}}
\newcommand{\ea}{\end{eqnarray}}
\newcommand{\vp}{\varphi}

\newcommand{\ve}{\varepsilon}

\begin{document} 
\begin{center}
{\Large{\sc On surfaces of prescribed weighted mean curvature}}\\[2cm]
{\large{\sc Matthias Bergner, Jens Dittrich}}\\[1cm]
{\small\bf Abstract}\\[0.4cm]
\begin{minipage}[c][2cm][l]{12cm}
{\small 
Utilizing a weight matrix 
we study surfaces of prescribed weighted mean curvature which
yield a natural generalisation to critical points 
of anisotropic surface energies. We first derive a differential
equation for the normal of immersions with prescribed weighted 
mean curvature, generalising a result of Clarenz and von der Mosel.
Next we study graphs of prescribed weighted mean curvature, for
which a quasilinear elliptic equation is proved. Using this equation,
we can show height and boundary gradient estimates. Finally, we solve
the Dirichlet problem for graphs of prescribed weighted mean curvature.}
\end{minipage}
\end{center}\noindent \\[2cm]
{\bf\large Introduction}\\[0.3cm]
Given some open set $U\subset\mathbb R^n$ let $X:U\to\mathbb R^{n+1}$
be a smooth immersion.
We denote by $N:U\to S^n$ its normal vector.
In \cite{clarenz2} Clarenz and von der Mosel studied 
critical points of the specific parametric functional
\bee\label{l1}
A(X)=\int_U \Big(F(N) + \langle Q(X),N\rangle \Big)dA \; .
\ee
Requiring the homogeneity condition $F(t p)=t F(p)$ for all $p\in\mathbb R^{n+1}$ and $t>0$,
this functional becomes invariant under reparametrisation of the surface.
The Euler equation of this functional leads
to surfaces $X$ whose weighted mean curvature $H_F$  
is prescribed by $H_F=\mbox{div} Q$.
A~simple example is the area functional with $F(p)=|p|$ and $Q\equiv 0$,
leading to surfaces whose mean curvature $H$ vanishes, i.e. minimal surfaces. 
In case of $F(p)=|p|$ together with some arbitrary $Q$ 
one obtains surfaces of prescribed mean curvature $H=\mbox{div} Q$. \\ \\
We will now generalise the class of prescribed weighted mean curvature
surfaces: We allow surfaces which do not necessarily arise as 
critical points of parametric functionals.
To this end, let us consider a symmetric $(n+1)\times (n+1)$ weight matrix 
$$ G=G(p):\mathbb R^{n+1}\backslash\{0\}\to\mathbb R^{(n+1)\times (n+1)} \; . $$
We require two conditions on the weight matrix:
First, the ellipticity condition
\bee\label{assumption_1}
\langle G(p) y,y\rangle >0 \quad \mbox{for all}\;\; y\in p^\perp\backslash\{0\}
\ee
i.e. $G(p)$ restricted to the $n$-dimensional space $p^\perp=\{y\in\mathbb R^{n+1}\, | \, \langle y,p\rangle=0\}$ 
is positive definite. Secondly we assume
\bee\label{assumption_2}
t G(t p)=G(p)\quad \mbox{and}\quad G(p)p=0 \quad \mbox{for all}\;\;p\in\mathbb R^{n+1}\backslash\{0\}
\; , \; t>0
\ee
i.e. $G(p)$ homogeneous of degree $-1$ and $p$ belongs to the kernel of $G(p)$.
Critical points of the functional (\ref{l1}) will be included 
in our considerations. For that case we just have to define the 
the weight matrix $G$ as the Hesse matrix of the second derivatives of $F$. Assumption (\ref{assumption_2})
on $G$ then follows directly from the $1$-homogeneity assumption on $F$. In particular,
the area functional $F(p)=|p|$ is included where the weight matrix is given by $G(p)=|p|^{-3}(E|p|^2-p p^T)$,
$E$ denoting the identity matrix. \\ \\
Similarly to Clarenz and von der Mosel (see \cite{clarenz1}, \cite{clarenz2})
we now define the weighted mean curvature of the surface $X$ as
$$ H_G:=\mbox{tr}(g^{-1}\, A_G\, g^{-1} b)=\mbox{tr}(g^{-1}\, A_G\, S) \; . $$
Here, the matrix $g$ is the first fundamental form defined by
$$ g:=DX^T DX \quad \mbox{with}\quad g_{ij}=\langle\partial_i X,\partial_j X\rangle $$
and the matrix $A_G$ is the {\it weighted first fundamental form}
$$ A_G:=DX^T G(N) DX \quad \mbox{with}\quad (A_G)_{ij}=\langle G(N)\, \partial_i X,\partial_j X\rangle \; . $$
Additionally, $b$ denotes the second fundamental form defined by
$$ b:=-DN^T DX \quad , \quad b_{ij}=-\langle \partial_i N,\partial_j X\rangle=\langle N,\partial_{ij} X\rangle \; . $$
Finally, $S:=g^{-1} b$ denotes the shape operator of the surface $X$. 
If we diagonalize $S$ at some fixed point on the surface, i.e. $S=\mbox{diag}(\kappa_1,\dots,\kappa_n)$
with the principal curvatures $\kappa_i$, then we obtain
\bee\label{l2}
H_G=\mbox{tr}(g^{-1}\, A_G\, S)=\sum_{i=1}^n \lambda_i(N) \kappa_i 
\ee
where $\lambda_i(N)$ are the diagonal entries of $g^{-1} A_G(N)$.
Hence, $H_G$ is a weighted sum of the principal curvatures of $X$. 
If the weight matrix $G$ is the identity on the tangent space, then
first fundamental form $g$ and weighted first fundamental form $A_G$ agree
and we obtain $\lambda_i(N)\equiv 1$ for $i=1,\dots,n$.
The weighted mean curvature then reduces to $H_G=\mbox{tr}(S)=\sum \kappa_i$, i.e. the
classical mean curvature of a surface. \\ \\
In this paper we study surfaces whose weighted mean curvature 
at every point $X$ is equal to some
prescribed function $\mathcal H=\mathcal H(X):\mathbb R^{n+1}\to\mathbb R$,
i.e. $H_G(X,N)=\mathcal H(X)$. As already mentioned, such surfaces arise for example 
as critical points of the functional
(\ref{l1}), where the prescribed weighted mean curvature is given by $\mathcal H(X):=\mbox{div} Q(X)$.
The special case $\mathcal H\equiv 0$ yields surfaces whose weighted mean curvature
vanishes. Such surfaces are called weighted minimal surfaces surfaces and yield
a natural extension of the class of minimal surfaces. Another interesting special case is
$\mathcal H\equiv \mbox{const}$, i.e. surfaces with constant weighted mean curvature as studied in \cite{froehlich},
a natural extension of cmc-surfaces. \\ \\
In Section 1 we start with a differential equation for the normal $N$ of
immersions of prescribed weighted mean curvature
(see Theorem \ref{t_normale}). It will then be used to derive a gradient maximum 
principle for graphs of prescribed weighted mean curvature (see Corollary \ref{c1}).
In Section 2 we derive a quasilinear, elliptic differential 
equation for graphs of prescribed weighted mean curvature.
This equation will be used to derive a height estimate for graphs,
using spherical caps as barriers.
In Section 3 we prove a boundary gradient estimate for graphs of prescribed weighted
mean curvature defined over a $C^2$-domain $\Omega\subset\mathbb R^n$. 
There, the weighted mean curvature of the boundary $\partial\Omega$ will play an important role.
Finally, we combine all these results to solve the Dirichlet problem for
graphs of prescribed weighted mean curvature in Section 4.
\subsection{The differential equation for the normal}
Given a parametrisation $X:\Omega\to\mathbb R^{n+1}$
with its normal $N:\Omega\to S^n$ we consider the first
and second fundamental forms defined by
$$ g_{ij}:=\langle \partial_i X,\partial_j X\rangle \quad , \quad
b_{ij}:=\langle\partial_{ij} X,N\rangle \quad i,j=1,\dots,n \; . $$
Let $g^{ij}$  be the inverse matrix of $g_{ij}$. 
We now define the matrix $A_G$ matrix with the entries
$$ a_{ij}:=\langle \partial_i X, G(N) \partial_j X\rangle $$
which is called the weighted first fundamental form.
Note $a_{ij}$ is a twice covariant tensor in the following sence: 
Under parameter transformations it transforms the same way as the first fundamental form $g_{ij}$.
If $G$ is the identity on the tangent space, then 
first and weighted first fundamental form agree. 
To derive a differential equation for the normal $N$ of the surface, we have to
consider the following differential operator.
\begin{definition}
Let $\psi\in C^2(\Omega,\mathbb R)$ be a function.
We define the weighted Laplace-Beltrami operator or $G$-Laplace-Betrami operator of $\psi$ by
$$ \triangle_G\psi:=\frac{1}{\sqrt{\det g}}\partial_i 
\Big(\sqrt{\det g} g^{ij} a_{jk} g^{kl} \partial_l \psi\Big) \; , $$
using the sum convention (summing from $1$ to $n$ over indices appearing twice).
\end{definition}
Remarks:
\begin{itemize}
\item[1.)] The weighted Laplace-Beltrami operator was already introduced
by Clarenz and von der Mosel in~\cite{clarenz1} and \cite{clarenz2}.
\item[2.)] In case that the weight matrix $G$ is the identity on the tangent space,
we have $a_{jk}=g_{jk}$ and the weighted Laplace-Betrami operator reduces
to the classical Laplace-Betrami operator.
\item[3.)] From the ellipticity assumption (\ref{assumption_1}) it follows that 
$a_{ij}$ is a positive definite matrix and so is $g^{ij} a_{jk} g^{kl}$.
Hence, $\triangle_G$ is an elliptic differential operator.
\end{itemize}
We now want to show a representation of the weighted Laplace-Betrami operator
in terms of the covariant derivative. For the derivation 
we will use the notations from Ricci calculus (see \cite[Chapter 4]{dubrovin}).
But let us first recall the definition
of covariant derivative.
We define the Christoffel symbols by
$\Gamma_{ij}^l:=\frac{1}{2}g^{lk}(\partial_i g_{jk}+\partial_j g_{ik}-\partial_k g_{ij})$.
Then the covariant derivative of a $1$-covariant tensor $T_i$ is defined by
\bee\label{l31}
D_i T_j:=\partial_i T_j-\Gamma^k_{ij} T_k \; , 
\ee
where $\partial_i$ denotes the usual derivative in direction $e_i$.
Secondly, we need the definition of the covariant derivative of a $2$-covariant, not necessarily symmetric tensor $T_{ij}$ by
\bee\label{l32}
D_k T_{ij}:=\partial_k T_{ij}-\Gamma_{ik}^l T_{lj}-\Gamma_{jk}^l T_{il} \; .
\ee
For the general definition of the covariant derivative of a tensor see 
\cite[Theorem 28.2.6]{dubrovin} or \cite[\S 63]{blaschke}.
The covariant derivative satisfies a product rule. 
Moreover the identity $D_i g_{jk}=D_i g^{jk}=0$ holds, known as the Lemma of Ricci.
Finally, in terms of the covariant derivative the Codazzi equations 
can be written as $D_i b_{jk}=D_k b_{ij}$.
We can now show
\begin{lemma}\label{lemma5}
The weighted Laplace-Betrami operator satisfies
$$ \triangle_G\psi= D_i\Big(g^{ij} a_{jk} g^{kl}\partial_l \psi\Big)
=g^{ij} \Big(D_i a_{jk}\Big) g^{kl} \partial_l\psi+g^{ij} a_{jk} g^{kl} D_{il} \psi $$
where $D_i$ denotes the covariant derivative of a tensor.
In particular, $\triangle_G \psi$ is parameter invariant.
\end{lemma}
{\it Proof:} \\
Let us set
$$ T^i:=g^{ij} a_{jk} g^{kl} \partial_l \psi $$
and note that $T^i$ is a $1$-contravariant tensor. 
In \cite[Example 29.3.4]{dubrovin} the following formula for the covariant derivative is proven
$$ D_i T^i=\frac{1}{\sqrt{\det g}}\partial_i \Big(\sqrt{\det g}\, T^i\Big) \; . $$
Using this, we obtain
$$ \triangle_G \psi=\frac{1}{\sqrt{\det g}} \partial_i\Big(\sqrt{\det g}\, T^i\Big)=D_i T^i \; .$$
This proves the first representation of $\triangle_G$ claimed in this lemma.
To prove the second one, we use the product rule as well as
the Lemma of Ricci to calculate
$$ \triangle_G\psi= D_i\Big(g^{ij} a_{jk} g^{kl}\partial_l \psi\Big)
=g^{ij} \Big(D_i a_{jk}\Big) g^{kl} \partial_l\psi+g^{ij} a_{jk} g^{kl} D_{il} \psi ,$$
ending the proof. \hfill $\Box$ \\ \\
We will now derive a differential equation for the normal vector $N$.
If $(N^1,\dots,N^{n+1})$ are the components of $N$  we define
$\triangle_G N:=(\triangle_G N^1,\dots,\triangle_G N^{n+1})$, i.e.
componentwise.
\begin{lemma}
Let $X\in C^3(\Omega,\mathbb R^{n+1})$ be a surface of prescribed weighted mean
curvature $\mathcal H\in C^1(\mathbb R^{n+1},\mathbb R)$.
Then its normal vector $N$ satisfies the differential equation
\bee\label{l22}
 \triangle_G N-g^{ij}(D_i a_{jk})g^{kl}\partial_l N
+\Big(\mbox{tr}(g^{-1} A_G S^2)-\langle \nabla \mathcal H,N\rangle\Big) N
=-\nabla\mathcal H+g^{ij}(D_p a_{jk}) g^{kl} b_{li} g^{pq}\partial_q X \; .
\ee
\end{lemma}\label{lemma6}
{\it Proof:} \\
From the assumption $X\in C^3(\Omega,\mathbb R^{n+1})$ we immediately
conclude the regularity $N\in C^2(\Omega,\mathbb R^{n+1})$.
We recall the Gauss-Weingarten equations
\bee\label{l23}
\partial_i N=-b_{ij} g^{jk}\partial_k X \quad \mbox{and} \quad
\partial_{ij} X=\Gamma_{ij}^k \partial_k X+b_{ij} N \; . 
\ee
Using the definition of the covariant derivative (\ref{l31}) of a $1$-covariant tensor we compute
$$ D_i \partial_j X=\partial_{ij} X-\Gamma_{ij}^k \partial_k X
=b_{ij} N+\Gamma_{ij}^k \partial_k X-\Gamma_{ij}^k \partial_k X=b_{ij} N \; . $$
Together with Lemma \ref{lemma5}, the product rule and the Lemma of Ricci we obtain
\bea
\triangle_G N
&=&g^{ij} (D_i a_{jk}) g^{kl} \partial_l N
+g^{ij} a_{jk} g^{kl} D_i \partial_l N \nonumber \\
&=&g^{ij} (D_i a_{jk}) g^{kl} \partial_l N
-g^{ij} a_{jk} g^{kl} D_i\Big(b_{lp} g^{pq} \partial_q X\Big) \nonumber \\
&=&g^{ij} (D_i a_{jk}) g^{kl} \partial_l N
-g^{ij} a_{jk} g^{kl} (D_i b_{lp}) g^{pq} \partial_q X 
-g^{ij} a_{jk} g^{kl} b_{lp} g^{pq} b_{qi} N \nonumber \\
&=&g^{ij} (D_i a_{jk}) g^{kl} \partial_l N
-g^{ij} a_{jk} g^{kl} (D_i b_{lp}) g^{pq} \partial_q X 
-\mbox{tr}(g^{-1} A_G S^2) N \label{l20}\; . 
\ea
For the second term in this sum we use the Codazzi equations $D_i b_{lp}=D_p b_{li}$ 
and the definition of weighted mean curvature 
$\mathcal H=g^{ij} a_{jk} g^{kl} b_{li}$ to get
\bea
g^{ij} a_{jk} g^{kl} (D_i b_{lp}) &=&
g^{ij} a_{jk} g^{kl} (D_p b_{li})
+g^{ij} (D_p a_{jk}) g^{kl} b_{li}
-g^{ij} (D_p a_{jk}) g^{kl} b_{li} \nonumber \\
&=&D_p\Big(g^{ij} a_{jk} g^{kl} b_{li}\Big)
-g^{ij} (D_p a_{jk}) g^{kl} b_{li} \nonumber \\
&=&D_p \mathcal H-g^{ij} (D_p a_{jk}) g^{kl} b_{li}
=\partial_p \mathcal H-g^{ij} (D_p a_{jk}) g^{kl} b_{li} \; . \label{l21}
\ea
Finally, using $\mathcal H=\mathcal H(X)$ and the chain rule we compute
$$ \partial_p \mathcal H g^{pq} \partial_q X=
\langle \nabla\mathcal H,\partial_p X\rangle g^{pq}\partial_q X
=\nabla\mathcal H-\langle \nabla\mathcal H,N\rangle N \; ,$$
noting that second term in this equation is the orthogonal projection
of $\nabla\mathcal H$ onto the tangent space.
Combining this with (\ref{l20}) and (\ref{l21}) yields the desired equation for
the normal $N$. \hfill $\Box$ \\ \\
{\it Remark:}
In case that the weight matrix is the identity on the tangent space
we obtain $a_{ij}=g_{ij}$ and by the Ricci Lemma $D_k a_{ij}=0$.
The differential equation then takes the form
$$ \triangle N+\Big(\mbox{tr}(S^2)-\langle\nabla\mathcal H,N\rangle\Big) N=-\nabla \mathcal H \; , $$
where now $\triangle$ is the classical Laplace-Betrami operator.
This is the well known differential equation for the normal vector
of a surface with (non-weighted) mean curvature $\mathcal H$. For dimension $n=2$ this equation
was proven by Sauvigny in \cite[Satz 1]{saudok} (see also \cite[Chapter XII, \S 9, Lemma 2]{saubuch3}). \\ \\
The differential equation (\ref{l22}) is not quite satisfying as it contains
a linear combination of the tangent vectors $\partial_1 X,\dots,\partial_n X$
on the right side. To get rid of this term, we will now replace it by a linear combination of
the derivatives $\partial_i N$, $i=1,\dots,n$ of the normal. 
If we assumed $K\neq 0$ for the Gaussian curvature of $X$, then we could
directly replace each $\partial_i X$ by a linear combination of $\partial_i N$,
as the vectors $\partial_i N$ would then be linearly independent. However,
we do not want to assume $K\neq 0$ as this is quite restrictive.
Instead, we will now use the character of the weight matrix $G=G(N)$
only depending on the normal $N$ but not on $X$ to achieve this substitution.
\begin{theorem}\label{t_normale}
Let $X\in C^3(\Omega,\mathbb R^{n+1})$ be an immersion of prescribed weighted mean
curvature $\mathcal H\in C^1(\mathbb R^{n+1},\mathbb R)$.
Then its normal vector $N$ satisfies the differential equation
$$ \triangle_G N+P^i \partial_i N+\Big(\mbox{tr}(g^{-1} A_G S^2)-\langle \nabla \mathcal H,N\rangle\Big) N
=-\nabla\mathcal H \quad \mbox{in}\; \;\Omega $$
with certain coefficients $P^i\in C^0(\Omega,\mathbb R)$. 
\end{theorem}
{\it Proof:}
1.) We first claim the following shift formula: 
For any vector $V\in\mathbb R^{n+1}$ we have the identity
\bee\label{l24}
\langle V,\partial_i N\rangle g^{ij} \partial_j X=\langle V,\partial_i X\rangle g^{ij} \partial_j N \; . 
\ee
To prove it, we use (\ref{l23}) and compute
$$ \langle V,\partial_i N\rangle g^{ij}\partial_j X=
-\langle V,\partial_l X\rangle b_{ik} g^{kl} g^{ij}\partial_j X  
=-\langle V,\partial_l X\rangle g^{kl} b_{ki} g^{ij}\partial_j X=\langle V,\partial_l X\rangle g^{lk}\partial_k N \; . $$
2.) Using the definition
$$ a_{ij}:=\langle \partial_i X, G(N)\partial_j X\rangle $$
together with (\ref{l23}) and the assumption $G(N)N=0$ we compute
\bea
\partial_k a_{ij}&=&\langle \partial_i X,\partial_k G(N) \partial_j X\rangle+
\langle \partial_{ik} X, G(N) \partial_j X\rangle
+\langle \partial_i X, G(N)\partial_{jk} X\rangle \nonumber \\
&=&\langle \partial_i X,\partial_k G(N)\partial_j X\rangle+\Gamma_{ik}^l a_{il}+\Gamma_{jk}^l a_{jk}\nonumber \; .
\ea
Together with the definition (\ref{l32}) of the covariant derivative of a twice covariant tensor we obtain
$$ D_k a_{ij}:=\partial_k a_{ij}-\Gamma_{ik}^l a_{li}-\Gamma_{jk}^l a_{il}
=\langle \partial_i X,\partial_k G(N)\partial_j X\rangle \; . $$
Now let $(N^1,\dots,N^{n+1})$ be the components of $N$. Then the chain rule gives
$$ D_k a_{ij}=\sum_{\mu=1}^{n+1} \langle \partial_i X, \partial_{N^\mu} G(N) \partial_j X\rangle \partial_k N^\mu 
=\langle V_{ij},\partial_k N\rangle \; ,$$
if we define the vectors $V_{ij}\in\mathbb R^{n+1}$, $i,j=1,\dots,n$, by
$$ V_{ij}^\mu:=\langle \partial_i X, \partial_{N^\mu} G(N) \partial_j X\rangle \quad \mbox{for}\; \;\mu=1,\dots,n+1 \; . $$
3.) Using 1.) and 2.) we can now rewrite the tangential term on the right side of (\ref{l22}) as
$$ g^{ij} (D_p a_{jk}) g^{kl} b_{li} g^{pq} \partial_q X
=g^{ij} g^{kl} b_{li} \langle V_{jk},\partial_p N\rangle g^{pq} \partial_q X
=g^{ij} g^{kl} b_{li} \langle V_{jk},\partial_p X\rangle g^{pq} \partial_q N \; . $$
If we define
$$ P^q:=-g^{ij} (D_i a_{jp}) g^{pq}-g^{ij} g^{kl} b_{li} \langle V_{jk},\partial_p X\rangle g^{pq} $$
then the theorem follows. \hfill $\Box$ \\ \\
{\it Remark:} If we consider the variational problem (\ref{l1}),
then weight matrix $G$ is obtained as the Hesse matrix
of some $C^3$-function $F:\mathbb R^{n+1}\backslash\{0\}\to\mathbb R$.
Using the Lemma of Schwarz, a computation reveals that
all $P^i$ vanish in that case. We then obtain the same differential equation as \cite[Corollary~4.3]{clarenz2}.
The proof in that paper relies on a formula for the second variation of the functional (\ref{l1}).
As our problems do not necessarily arise as Euler equations of 
variational problems, we do not have the tool of second variation at hand. 
Instead, we have used only geometric identities to derive our equation for the normal. \\ \\
We now use the differential equation to derive a gradient maximum principle
for graphs of prescribed weighted mean curvature.
We need the following inequality
\bee\label{l25}
\mbox{tr}(g^{-1} A_G S^2) \mbox{tr}(G)\geq \mbox{tr}(g^{-1} A_G S)^2=(H_G)^2 \; .
\ee
As this inequality is invariant under repametrisation, it suffices to prove
it for a particular parametrisation. Given some point $p_0\in\mathbb R^{n+1}$ on the surface, let
$X:B\to\mathbb R^{n+1}$ be a parametrisation satisfying
$X(0)=p_0$, $g_{ij}(0)=\delta_{ij}$ and $b_{ij}(0)=S^i_j(0)=\mbox{diag}(\kappa_1,\dots,\kappa_n)$
with the principal curvatures $\kappa_i$ of the surface at $p_0$.
Using Cauchy-Schwarz inequality we can then estimate
$$ \mbox{tr}(g^{-1} A_G S)^2=\Big(\sum_{i=1}^n a_{ii} \kappa_i\Big)^2
=\Big(\sum_{i=1}^n \sqrt{a_{ii}}\, \sqrt{a_{ii}}\kappa_i\Big)^2
\leq \Big(\sum_{i=1}^n a_{ii}\Big)\Big(\sum_{i=1}^n a_{ii}\kappa_i^2\Big)
=\mbox{tr}(g^{-1} A_G S^2)\,\mbox{tr}(G) \; . $$
Here we use $a_{ii}=\langle \partial_i X, G(N)\partial_i X\rangle\geq 0$
which follows from the ellipticity assumption (\ref{assumption_1}) on $G$.
Secondly, we use $\sum a_{ii}=\mbox{tr}(G)$ (see the proof of Lemma \ref{lemma2}). 
We can now prove
\begin{corollary}\label{c1}
For $u\in C^3(\Omega,\mathbb R)\cap C^1(\ol\Omega,\mathbb R)$
let $X(x)=(x,u(x))$ be a graph of prescribed weighted mean curvature
$\mathcal H=\mathcal H(x,z)\in C^1(\ol\Omega\times\mathbb R,\mathbb R)$
satisfying the monotonocity assumption
$\frac{\partial}{\partial z} \mathcal H\geq 0$.
Additionally, we require
\bee
\mathcal H^2(x,z)\geq \mbox{tr}(G(p))|\nabla\mathcal H(x,z)| \quad \mbox{for all}\; \;
x\in\Omega\; , \; z\in\mathbb R\; , \; p\in S^n \; . \label{l26}
\ee
Then the estimate holds
$$ \sup\limits_{\Omega} |\nabla u|\leq \sup\limits_{\partial\Omega} |\nabla u| \; . $$
\end{corollary}
{\it Proof:} \\
Consider the last component of the normal $\psi(x):=N^{n+1}(x)=(1+|\nabla u|^2)^{-1/2}>0$.
By Theorem~\ref{t_normale} together with the assumption $\mathcal H_z\geq 0$ it satisfies
the differential inequality
$$ \triangle_G \psi+P^i(x)\partial_i \psi+\Big(\mbox{tr}(g^{-1} A_G S^2)
-|\nabla\mathcal H|\Big)\psi \leq 0 \; . $$
Using (\ref{l25}) together with assumption (\ref{l26}) then yields
$$ \triangle_G \psi+P^i(x) \partial_i \psi\leq 0 \quad \mbox{in}\; \Omega \; .  $$
By the maximum principle $\psi$ achieves its minimum on $\partial\Omega$
and hence $|\nabla u|$ must achieve its maximum on $\partial\Omega$. \hfill $\Box$ \\ \\
{\it Remark:} We have used the differential equation for the normal to derive
a gradient maximum principle. Aside from this, the it
may also be for other things. For example, it may be used to derive
purely interior a priori gradient estimates for graphs.
Within the context of the functional (\ref{l1}),
the differential equation for the normal is used in \cite[Theorem 1.4]{clarenz2}
to prove a projectability theorem. This result states that under certain geometric conditions
any stable, immersed parametric surface of prescribed weighted mean curvature
must be a graph over the $x_1,x_2$-plane.
\subsection{Graph representation and $C^0$-estimate}
For a function $u:\Omega\to\mathbb R$ let us consider
the graph parametrisation $X(x):=(x,u(x))$ together with the
upper normal vector
$$ N(x):=\frac{1}{\sqrt{1+|\nabla u|^2}} (-\nabla u,1) \quad \mbox{for}\;\; x\in\Omega \; . $$
Then we say that $u$ is a graph of prescribed weighted mean curvature if its
parametrisation $X(x):=(x,u(x))$ has prescribed weighted mean curvature.
\begin{lemma}\label{lemma1}
Let $u\in C^2(\Omega,\mathbb R)$ be a graph of prescribed weighted mean
curvature $\mathcal H:\Omega\times\mathbb R\to\mathbb R$. Then $u$ satisfies the quasilinear, elliptic
differential equation 
\bee\label{l74}
\sum_{i,j=1}^n G_{ij}(-\nabla u,1) \partial_{ij} u=\mathcal H(x,u) \quad \mbox{in}\;\; \Omega \; . 
\ee
\end{lemma}
{\it Proof:} \\
For the parametrisation $X(x):=(x,u(x))$ the first fundamental form
$g=g_{ij}$ is given by
$$ g=E+\nabla u\nabla u^T \quad , \quad g_{ij}=\langle \partial_i X,\partial_j X\rangle
=\delta_{ij}+\partial_i u\partial_j u \; .$$
Next we compute the second fundamental form $b=b_{ij}$ as
$$ b_{ij}=\langle \partial_{ij} X, N\rangle=\frac{\partial_{ij} u}
{\sqrt{1+|\nabla u|^2}} \; . $$
Now let 
$$ a_{ij}:=\langle G(N)\, \partial_i X,\partial_j X\rangle $$
be the entries of the matrix $A_G=(DX)^T G(N) DX$.
From $\partial_i X=(e_i,\partial_i u)$ we obtain the representation
\bee\label{l8}
a_{ij}=G_{ij}+G_{i\, n+1}\partial_j u+G_{j\, n+1}\partial_i u+G_{n+1\,n+1}\partial_i u\partial_j u \quad
\mbox{for}\; \;i,j=1,\dots,n \; .
\ee
Let us now decompose the matrix $G$ into 
\bee
G = \left (\begin{array}{cc}\hat G& w\\w^T& c \end{array}\right )
\ee
where $\hat G$ are the first $n\times n$ entries of $G$ and $(w^T,c)$ is the last row of $G$.
Noting that $(\nabla u,-1)$ is a multiple of the normal $N$, assumption
(\ref{assumption_2}), i.e. $G(N)N=0$, leads to
$$\hat G \nabla u=w \quad \mbox{and} \quad \nabla u^T w=c= w^T\nabla u \; . $$
Using this we compute
\bea
 g\,\hat G\,g &=& (E+\nabla u\nabla u^T)\,\hat G\, (E+ \nabla u\nabla u^T)\nonumber\\
 &=& \hat G + \nabla u w^T + w \nabla u^T + c\nabla u\nabla u^T\stackrel{(\ref{l8})}{=} A_G\nonumber \; ,
\ea
which is equivalent to $\hat G=g^{-1}\, A_G\, g^{-1}$.
Employing $S=g^{-1} b$ for the shape operator together with the definition of the weighted mean curvature we then obtain
\bea
H_G&=&\mbox{tr}(g^{-1}\, A_G S)=\mbox{tr}(g^{-1}\, A_G\, g^{-1} b)
 =\mbox{tr}(\hat G\, b)=\sum_{i,j=1}^n G_{ij}(N)\frac{\partial_{ij} u}{\sqrt{1+|\nabla u|^2}}\nonumber \\
&=&\sum_{i,j=1}^n G_{ij}(-\nabla u,1) \partial_{ij} u \nonumber \; .
\ea
In the last step we have used  the $-1$-homogeneity assumption (\ref{assumption_2})
on the weight matrix $G$. \hfill $\Box$ \\ \\
{\it Remarks:}
\begin{itemize}
\item[1.)] Note that only the first $n\times n$
entries $G_{ij}$ for $i,j=1,\dots,n$ enter into the differential equation.
This is due to the symmetry assumption $G^T=G$ and assumption (\ref{assumption_2}) $G(p)p=0$.
Indeed, once the first $n\times n$ entries of $G$ are given, the remaining entries
$G_{i\, n+1}$ and $G_{n+1\, i}$ are uniquely determined by the above relations.
\item[2.)] Assuming $G=G(p)$ to be a differentiable 
function of $p$, the maximum and comparision principle of \cite[Theorem 10.1]{gilbarg}
can be applied to solutions $u$ of equation (\ref{l74}).
\item[3.)] For the special choice of the weight
matrix $G(p)=|p|^{-3}(E|p|^2-p p^T)$, corresponding to functional (\ref{l1}) with $F(p)=|p|$,
the differential equation takes the form
$$ (1+|\nabla u|^2)^{-3/2}
\Big ((1+|\nabla u|^2)\delta_{ij}-\partial_i u\partial_j u\Big)\partial_{ij} u
= \mathcal H(x,u) \quad \mbox{in}\; \; \Omega \; . $$
This is the classical nonparametric equation for a graph of prescribed (non-weighted) 
mean curvature $\mathcal H(x,u)$. 
\item[4.)] As we can see from the example above, the quasilinear elliptic
equation under consideration is not uniformly elliptic!
\end{itemize}
The next example will illustrate that graphs of prescribed weighted mean curvature
are obtained as critical points of certain geometric, nonparametric functionals. 
\begin{example}\label{example1}
Given two functions $F\in C^2(\mathbb R^{n+1}\backslash\{0\},\mathbb R)$
and $b\in C^1(\ol\Omega\times\mathbb R,\mathbb R)$, 
consider the nonparametric version of the functional (\ref{l1})
\bee\label{l70}
A(u):=\int_\Omega\Big(F(-\nabla u,1)+b(x,u)\Big) dx 
\ee
whose Euler equation is given by
$$ \sum_{i,j=1}^n F_{p_i p_j}(-\nabla u,1) \partial_{ij} u =b_z(x,u) \quad \mbox{in}\;\; \Omega \; . $$
This is exactly the differential equation of Lemma \ref{lemma1}
if we define 
$$ \mathcal H(x,z):=b_z(x,z) \quad \mbox{and}\quad G_{ij}(p):= F_{p_i p_j}(p) \quad \mbox{for}\;\; i,j=1,\dots,n+1 \; . $$
Hence, critical points of the functional (\ref{l70}) can be interpreted as
graphs with prescribed weighted mean curvature.
This weight matrix $G$ will satisfy both of the required assumptions (\ref{assumption_2})
if we assume $F$ to be $1$-homogeneous, i.e. $F(tp)=t F(p)$ for all $t>0$.
\end{example}
\begin{example}
A particularly interesting example is
$F(p)=\sqrt{p_1^2+\dots+p_n^2+\ve^2 p_{n+1}^2}$ with the corresponding functional 
$$ A_\ve(u)=\int_\Omega\Big( \sqrt{\ve^2+|\nabla u|^2}+b(x,u)\Big) dx $$
for $\ve>0$. This functional can be viewed as a regularised version of the functional
$$ A_0(u)=\int_\Omega\Big(|\nabla u|+b(x,u)\Big) dx \; . $$
The Euler equation of this functional $A_0$ is the degenerated elliptic equation
$$ \mbox{div}\frac{\nabla u}{|\nabla u|}=b_z(x,u) \quad \mbox{in}\; \;\Omega $$
which is only welldefined if $\nabla u\neq 0$. A solution $u$ has the property that
its level sets $M_c:=\{x\in\Omega \, | \, u(x)=c\}\subset\mathbb R^n$ have prescribed
mean curvature $b_z(x,c)$ for $x\in M_c$. Hence, one obtains a family 
of implicitely defined surfaces having prescribed mean curvature.
\end{example}
We now want to derive an estimate of $C^0$-norm
for graphs using spherical caps as upper and lower barriers.
To do this, we first have to compute the weighted mean curvature of a sphere.
\begin{lemma}\label{lemma2}
Let $S_R=\{y\in\mathbb R^{n+1}\, : \, |y-y_0|=R\}$ be a sphere
of radius $R>0$. Then its weighted mean curvature $H_G(y)$ 
at some point $y\in S_R$ is given by $H_G(y)=\frac{1}{R}\, \mbox{tr}\, G(N)$, if $N$ is the inner
normal to $S_R$ at $y$ and by $H_G(y)=-\frac{1}{R}\, \mbox{tr}\, G(N)$, if $N$ is the outer normal to $S_R$ .
\end{lemma}
{\it Proof:} \\
Note that $S=\pm \frac{1}{R} E$ for the shape operator of the sphere $S_R$,
the sign depending on the choice of normal. Hence, we compute the weighted mean curvature by
$$ H_G(y)=\pm \mbox{tr}(g^{-1}\, A_G\, S)
=\pm \frac{1}{R}\mbox{tr}(g^{-1}\, A_G) \; . $$
The lemma now follows if $\mbox{tr}(g^{-1} A_G)=\mbox{tr}\, G$ holds.
To show this, let $X:B\to S_R$ be a parametrisation of $S_R$ 
with $X(0)=y$. Additionally we may assume that $g_{ij}(0)=\delta_{ij}=g^{ij}(0)$.
At the point $X(0)=y$ we then obtain
$$ \mbox{tr}(g^{-1} A_G)=\mbox{tr}\, A_G=
\langle G(N) \partial_i X,\partial_i X\rangle=
\langle G(N) \partial_i X,\partial_i X\rangle+\langle G(N) N, N\rangle
=\mbox{tr} G(N) \; . $$
Here we have used the assumption $G(p)p=0$ together with the fact that
$\partial_1 X(0),\dots,\partial_n X(0),N(0)$ form an orthonormal basis of $\mathbb R^{n+1}$,
which follows directly from $g_{ij}(0)=\delta_{ij}$. \hfill $\Box$ \\
\begin{theorem}($C^0$-estimate) \label{t1} \\
Let $u\in C^2(\Omega,\mathbb R)\cap C^0(\ol\Omega,\mathbb R)$
be a graph of prescribed weighted mean curvature $\mathcal H:\ol\Omega\times\mathbb R\to\mathbb R$
over a bounded domain $\Omega\subset B_R(0)=\{x\in\mathbb R^n\, : \, |x|<R\}$.
We assume the smallness condition
\bee\label{l3} 
|\mathcal H(x,z)|\leq \frac{1}{R}
\mbox{tr}\, G(p) \quad \mbox{for all}\; \; x\in\Omega \; , \; z\in\mathbb R\; , \; p\in S^n \; .
\ee
Then the following estimate holds
$$ \sup_\Omega |u(x)|\leq \sup_{\partial\Omega} |u(x)|+R \; . $$
\end{theorem}
{\it Proof:} \\
Let us define a spherical cap of radius $R$ by
$$ v(x):=-\sup\limits_{x'\in\partial\Omega} |u(x')|-\sqrt{R^2-|x|^2} \quad \mbox{for}\;\; x\in\ol\Omega $$
which is well definded because of $\Omega\subset B_R(0)$. 
Then Lemma \ref{lemma1}, Lemma \ref{lemma2} together with the smallness assumption (\ref{l3})
yield the differential inequality
$$ \sum_{i,j=1}^n G_{ij}(-\nabla v,1) \partial_{ij} v\geq \sum_{i,j=1}^n
G_{ij}(-\nabla u,1)\partial_{ij} u \quad \mbox{in}\;\; \Omega \; . $$
Noting $u\geq v$ on $\partial\Omega$, the comparision principle for 
quasilinear elliptic equations \cite[Theorem 10.1]{gilbarg} yields $u\geq v$ in $\Omega$.
Similarly we can show $-u\leq -v$ in $\Omega$, which then yields
$$ \sup_\Omega |u(x)|\leq \sup_\Omega |v(x)|
\leq \sup_{\partial\Omega} |u(x)|+R \; , $$
proving the claimed estimate. \hfill $\Box$ \\ \\
We can also use the maximum principle to prove a non-existence theorem.
\begin{theorem}(non-existence of graphs) \label{nonexistence}\\
Let $\Omega:=B_R(0)$ be the ball of radius $R>0$ centered at $0$.
Let a prescribed weighted mean curvature
$\mathcal H\in C^0(\Omega\times\mathbb R,\mathbb R)$ be given such that
\bee\label{l73}
\mathcal H(x,z)>\frac{1}{R} \mbox{tr}\, G(p)\geq 0 
\quad \mbox{for all}\;\; x\in\Omega \, , \, z\in\mathbb R \, , \, p\in S^n \; . 
\ee
Then a graph $u\in C^2(\Omega,\mathbb R)\cap C^1(\ol\Omega,\mathbb R)$ of prescribed
weighted mean curvature $\mathcal H$ does not exist.
\end{theorem}
{\it Proof:} \\
Assume to the contrary that such a graph $u\in C^2(\Omega,\mathbb R)\cap C^1(\ol\Omega,\mathbb R)$ exists. Let us now define
$$ v(x):=c-\sqrt{R^2-|x|^2} \quad \mbox{for}\;\; x\in \ol\Omega \; , $$
where $c$ is the smallest real number for which $u(x)\leq v(x)$ in $\Omega$.
Then there exists some $x_*\in \ol\Omega$ with $u(x_*)=v(x_*)$.
Lemma \ref{lemma1}, Lemma \ref{lemma2} together with the assumption (\ref{l73})
yield the differential inequality
$$ \sum_{i,j=1}^n G_{ij}(-\nabla v,1) \partial_{ij} v<\sum_{i,j=1}^n
G_{ij}(-\nabla u,1)\partial_{ij} u \quad \mbox{in}\; \;\Omega \; . $$
The comparision principle \cite[Theorem 10.1]{gilbarg}
then implies $x_*\in\partial\Omega$, i.e. $|x_*|=R$.
On the other hand, $u(x)\leq v(x)$ in $\Omega$ and $u(x_*)=v(x_*)$
imply $\frac{\partial u}{\partial \nu}(x_*)\geq \frac{\partial v}{\partial \nu}(x_*)$,
where $\nu$ is the outer normal to $\partial\Omega$ at $x_*$.
However, we have $\frac{\partial v}{\partial\nu}(x_*)=+\infty$ because of $|x_*|=R$,
contradicting $u\in C^1(\ol\Omega,\mathbb R)$.\hfill $\Box$ \\ \\
Note that this result can easily be generalised to domains $\Omega\subset\mathbb R^n$
satisfying $B_R(0)\subset\Omega$.
\subsection{Boundary gradient estimate for graphs}
In this section we will derive a boundary gradient estimate for
graphs of prescribed weighted mean curvature.
Roughly speaking, we will use the cylinder 
$Z_{\partial\Omega}=\{(x,z)\in\mathbb R^{n+1}\, : \, x\in\partial\Omega\}$
as the barrier. We will have to require that this cylinder has a sufficienlty 
large weighted mean curvature w.r.t the inner normal and sufficiently small
weighted mean curvature w.r.t. the outer normal.
A technical difficulty arises 
from the fact that the cylinder is not a graph over the $x_1,\dots,x_n$ hyperplane.
Instead, we will use as barrier a graph which is sufficiently close to the cylinder. 
This barrier will be defined in terms of the distance function
$d(x)=\mbox{dist}(x,\partial\Omega)$. For a $C^2$-domain $\Omega$
the distance function be of class $C^2$ within the set
$\Omega_\mu:=\{x\in\ol\Omega\, |\, d(x)<\mu\}$ for sufficiently small
$\mu=\mu(\Omega)>0$ (see \cite[Lemma 14.16]{gilbarg}).
To start, we first show
a formula which expresses the weighted mean curvature of the boundary
in terms of the distance function.
\begin{lemma}\label{lemma3}
Let $\Theta\subset\mathbb R^{n+1}$ be a $C^2$-domain,
$d(y):=\mbox{dist}(y,\partial\Omega)$ be the distance function.
Then the weighted mean curvature $H_G^+$ of $\partial\Theta$ w.r.t. the inner normal is given by
$$ H_G^+(y,\partial\Theta)=-\sum_{i,j=1}^{n+1} G_{ij}(\nabla d(y)) \partial_{ij} d(y) 
\quad \mbox{for}\;\; y\in\partial\Theta $$
while the weighted mean curvature $H_G^-$ of $\partial\Theta$ w.r.t. the outer normal is given by
$$ H_G^-(y,\partial\Theta)=\sum_{i,j=1}^{n+1} G_{ij}(-\nabla d(y)) \partial_{ij} d(y) 
\quad \mbox{for}\;\; y\in\partial\Theta \; . $$
\end{lemma}
{\it Proof:} \\
1.) To give the proof of the lemma, we first have to recall some
facts about the distance function $d(y)$. 
At first we have $|\nabla d|^2\equiv 1$.
Differentiating this equation yields 
\bee\label{l10}
\sum_{i=1}^{n+1} \partial_i d\, \partial_{ij} d=0 \quad \mbox{for all}\;\;j=1,\dots,n+1 \; . 
\ee
Next, for all $y\in\partial\Theta$ the gradient $\nabla d(y)$ is
equal to the interior unit normal to $\partial\Theta$. \\ \\
2.) To prove the lemma at some point $y_0\in\partial\Theta$,
we can assume $y_0=0$ after a suitable translation.
After an additional rotation in $\mathbb R^{n+1}$ we may locally represent
$\partial\Theta$ as a graph in the form $(x,\psi(x))$
with $0=\psi(0)=\nabla \psi(0)$. Additionally, we assume that the interior
normal to $\partial\Theta$ at $0$ is the vector $e_{n+1}:=(0,\dots,0,1)$, i.e. $\Theta$
lies above the graph $(x,\psi(x))$. By 1.) we also have $\nabla d(0)=e_{n+1}$, in particular
$\partial_i d(0)=0$ for all $i=1,\dots,n$. Putting this into 
(\ref{l10}) we obtain 
$$ \partial_{i\,n+1} d(0)=\partial_{n+1\,i} d(0)=0 \quad \mbox{for all}\;\; i=1,\dots,n+1 \; . $$
Now by twice differentiating the identity 
$d(x,\psi(x))=0$ and evaluating $0$ one gets $\partial_{ij} d(0)=-\partial_{ij} \psi(0)$.
Using Lemma \ref{lemma1} we compute the weighted mean curvature
of $\partial\Theta$ at $0$ w.r.t. the inner normal by
$$ H_G(0)=\sum_{i,j=1}^n G_{ij}(e_{n+1})\partial_{ij} \psi(0)
=-\sum_{i,j=1}^n G_{ij}(e_{n+1}) \partial_{ij} d(0)
=-\sum_{i,j=1}^{n+1} G_{ij}(e_{n+1}) \partial_{ij} d(0) \; , $$
proving the formula for $H_G^+$. Similarly we can derive the formula for
the weighted mean curvature w.r.t. the outer normal. \hfill $\Box$ \\ \\
{\it Remark:} In general, one cannot expect any kind of relation
between the quantities $H_G^+$ and $H_G^-$. However,
if we require the condition $G(-p)=G(p)$, the two weighted mean curvatures 
satisfy $H^+_G=-H^-_G$. This condition holds for example
in case of the usual (non-weighted) mean curvature 
where the weight matrix is given by $G(p)=|p|^{-3}(|p|^2 E-p p^T)$. 
\begin{definition}\label{def1}
Let $\Omega\subset\mathbb R^n$ be a $C^2$-domain.
Then we define the weighted mean curvature of $\partial\Omega$
w.r.t. the inward (or outward) normal to be the weighted mean curvature
of the boundary $\partial Z$ of the cylinder $Z_\Omega:=\{(x,z)\in\mathbb R^{n+1}\, : \, x\in\Omega\}$
w.r.t. to the inward (or outward) normal. 
In terms of the distance function $d(x)=\mbox{dist}(x,\partial\Omega)$,
they can be computed by
$$ H_G^+(x,\partial\Omega)=-\sum_{i,j=1}^n G_{ij}(\nabla d(x),0) \partial_{ij} d(x) 
\quad \mbox{for}\;\; x\in\partial\Omega $$
as well as
$$ H_G^-(x,\partial\Omega)=\sum_{i,j=1}^n G_{ij}(-\nabla d(x),0) \partial_{ij} d(x)  
\quad \mbox{for}\; \;x\in\partial\Omega \; .$$
\end{definition}
Now let $u\in C^2(\ol\Omega,\mathbb R)$ be a graph of prescribed weighted mean
curvature having Dirichlet
boundary values $u=\vp$ on $\partial\Omega$ for some $\vp\in C^2(\ol\Omega,\mathbb R)$.
To obtain upper and lower barriers for $u$, let us define
$$ v(x):=c\, d(x)+\vp(x) \quad \mbox{for}\; x\in \Omega_\mu:=\{x\in\ol\Omega\, | \, d(x)\leq \mu\} $$
for some constant $c\in\mathbb R$ and $\mu>0$. 
Fow sufficiently large $c>0$ we will obtain a upper and for sufficiently
small $c<0$ a lower barrier.
The upper unit normal of the graph $v$ is given by
$$ N_v:=\frac{1}{\sqrt{1+|\nabla v|^2}}(-\nabla v,1)=
\frac{1}{\sqrt{1+c^2+2c\langle \nabla d,\nabla \vp\rangle+|\nabla \vp|^2}} (-c \, \nabla d-\nabla \vp,1) \; . $$
Note the following convergence of the normal
$$ \nu(x):=\lim\limits_{c\to\infty} N_v(x)=(-\nabla d(x),0) \quad \mbox{in}\; \;\Omega_\mu\; . $$
For $x\in\partial\Omega$ the limit $\nu(x)$ is actually 
equal to the outer unit normal to $\partial\Omega$ at $x$. 
Using Lemma \ref{lemma1} we compute the weighted mean curvature of the graph $v$ by
$$ H_G(v)=G_{ij}(-\nabla v,1) \partial_{ij} v
=\frac{1}{\sqrt{1+|\nabla v|^2}} G_{ij}(N_v)\partial_{ij} v \; . $$
We can now compute the limit 
\bea
\lim_{c\to+\infty} H_G(v)(x)&=&\lim_{c\to+\infty} 
\frac{1}{\sqrt{
1+c^2+2\langle \nabla d,\nabla \vp\rangle+|\nabla \vp|^2}} G_{ij}(N_v)
(c\,\partial_{ij} d+\partial_{ij}\vp) \nonumber \\
&=&\lim_{c\to+\infty} \frac{c}{\sqrt{1+c^2+2\langle \nabla d,\nabla \vp\rangle+|\nabla \vp|^2}} G_{ij}(N_v)
\partial_{ij}d \nonumber \\
&=&G_{ij}(\nu) \partial_{ij} d=H_G^-(x) \quad \mbox{for}\;\; x\in\partial\Omega \;  \label{l12}
\ea
where $H_G^-$ is the weighted mean curvature of $\partial\Omega$ w.r.t. the outward normal
(see Definition \ref{def1}).
Similarly, for $c\to-\infty$ one gets the limit
$$ \lim\limits_{c\to-\infty} H_G(v)(x)=H_G^+(x) \quad \mbox{for}\; \;x\in\partial\Omega $$
with the weighted mean curvature $H_G^+$ w.r.t. the inner normal.
Combining these results we can show
\begin{theorem}\label{t_boundarygradient}(Boundary gradient estimate) \\
Assumptions:
\begin{itemize}
\item[a)] For some $C^2$-domain $\Omega\subset\mathbb R^n$ let $u\in C^2(\ol\Omega,\mathbb R)$
be a graph of prescribed weighted mean curvature $\mathcal H\in C^0(\ol\Omega\times\mathbb R,\mathbb R)$.
\item[b)] Assume that $u$ satisfies the boundary condition $u=\vp$ on $\partial\Omega$
for some $\vp\in C^2(\partial\Omega,\mathbb R)$. Additionally, we require 
the estimate $|u(x)|\leq M$ in $\Omega$ with some constant $M$.
\item[c)] Let $H_G^+:\partial\Omega\to(0,+\infty)$ be the weighted mean curvature
of $\partial\Omega$ w.r.t the inner normal and 
$H_G^-:\partial\Omega\to(-\infty,0)$ the weighted mean curvature w.r.t. the outer normal.
We then require
\bee\label{l11}
H_G^-(x)< \mathcal H(x,z)< H_G^+(x) \quad \mbox{for all}\;\; x\in\partial\Omega \, , \, |z|\leq M \; . 
\ee
\end{itemize}
Then we have the estimate
$$ \sup_{x\in\partial\Omega} |\nabla u(x)|\leq C $$
with a constant $C$ only depending on the data $\Omega$, $||\vp||_{C^2(\partial\Omega)}$, $M$ and 
the moduli of continuity of $\mathcal H$ and $G$.
\end{theorem}
{\it Proof:}
Given $\vp\in C^2(\partial\Omega,\mathbb R)$ we can extend it to
$\tilde\vp\in C^2(\ol\Omega,\mathbb R)$ (see \cite[Lemma 6.37]{gilbarg})
such that $\tilde\vp=\vp$ on $\partial\Omega$.
As above, let us consider $v(x)=c\, d(x)+\tilde \vp(x)$ for
$x\in\Omega_\mu$. Note that $v(x)=\vp(x)=u(x)$ on $\partial\Omega$.
By assumption (\ref{l11}) together with the limit (\ref{l12})
we can first determine $\mu>0$ and $c_0>0$ such that
$$ H_G(v)(x)<\mathcal H(x,z) \quad \mbox{for all}\;\; x\in\Omega_\mu \; 
, \; |z|\leq M \; , \; c\geq c_0 \; . $$
In particular, we have
$$ H_G(v)(x)<\mathcal H(x,u(x)) \quad \mbox{for}\;\; x\in\Omega_\mu \; , \; c\geq c_0 \; . $$
Defining $c_1:=(M+||\tilde\vp||_{C^0(\Omega)})\mu^{-1}$  we 
obtain $v(x)\geq M$ whenever $c\geq c_1$ and $d(x)=\mu$.
In particular, this implies $v(x)\geq u(x)$ on $\partial\Omega_\mu$
whenever $c\geq c_1$. Setting $c:=\max(c_0,c_1)$, the comparision 
principle for quasilinear elliptic equations yields
$v(x)\geq u(x)$ in $\Omega_\mu$ and
$$ \frac{\partial v(x)}{\partial \nu}\leq \frac{\partial u(x)}{\partial \nu}
\quad \mbox{on} \;\; \partial\Omega \; , $$
where $\nu$ is the outer unit normal to $\partial\Omega$ at $x$.
Similarly, we can construct a lower barrier by choosing $c<0$ sufficiently small.
This will yield an estimate of $|\frac{\partial u(x)}{\partial\nu}|$
and together with the Dirichlet boundary condition $u=\vp$ on $\partial\Omega$
we can give an estimate of $|\nabla u(x)|$ on $\partial\Omega$. \hfill $\Box$ \\ \\
{\it Remark:}
The methods we use are quite similar to \cite[Chapter 14.3]{gilbarg},
where boundary gradient estimates for general quasilinear elliptic equations
under boundary curvature conditions are derived. 
In fact, we could conclude the boundary gradient estimate also from
\cite[Theorem 14.9]{gilbarg}. There, certain structure conditions
on the differential operator are required. Those structure conditions
can be verified to hold for our problem using the homogeneity assumption 
$tG(tp)=G(p)$ on the weight matrix $G$.
\subsection{The Dirichlet problem for graphs of prescribed weighted mean curvature}
In this section we study the Dirichlet problem for graphs of prescribed weighted mean 
curvature: Given a weight matrix 
$$ G\in C^{1+\alpha}(\mathbb R^{n+1}\backslash\{0\},\mathbb R^{(n+1)\times (n+1)}) $$
satisfying the assumptions
(\ref{assumption_1}) and (\ref{assumption_2}), a prescribed weighted mean curvature
$\mathcal H\in C^{1+\alpha}(\ol\Omega\times\mathbb R,\mathbb R)$ and Dirichlet boundary values
$g\in C^{2+\alpha}(\partial\Omega,\mathbb R)$ we look for a solution of
\bee\label{dirichletproblem}
 u\in C^{2+\alpha}(\ol\Omega,\mathbb R) \quad , \quad
\sum_{i,j=1}^n G_{ij}(-\nabla u,1) \partial_{ij} u=\mathcal H(x,u) \quad \mbox{in}\;\; \Omega 
\quad \mbox{and} \quad  u=g \quad \mbox{on}\;\; \partial\Omega \; .
\ee
Combining the results we have proven so far we obtain
\begin{theorem}
Assumptions:
\begin{itemize}
\item[a)] Let $\Omega\subset\mathbb R^n$ be a $C^{2+\alpha}$-domain such that $\Omega\subset B_R(0)$ for some $R>0$.
Let $H_G^+:\partial\Omega\to(0,+\infty)$ be the weighted mean curvature
of $\partial\Omega$ w.r.t the inner normal and 
$H_G^-:\partial\Omega\to(-\infty,0)$ be the weighted mean curvature w.r.t. the outer normal.
\item[b)] Let $\mathcal H=\mathcal H(x,z)\in C^{1+\alpha}(\ol\Omega\times\mathbb R,\mathbb R)$
be the prescribed weighted mean curvature satisfying the monotonocity assumption $\mathcal H_z\geq 0$.
\item[c)] Let the inequalities 
\bea
&& R |\mathcal H(x,z)|\leq \mbox{tr}(G(p)) \quad , \label{a1} \\
&&H_G^-(x)<\mathcal H(x,z)<H_G^+(x) \quad \mbox{and} \label{a2} \\
&& \mathcal H^2(x,z)\geq \mbox{tr}(G(p)) |\nabla\mathcal H(x,z)| \quad \mbox{for all}\;\; x\in\ol\Omega \; , \; 
z\in\mathbb R \; , \; p\in S^n  \label{a3}
\ea
be satisfied.
\end{itemize}
Then for all Dirichlet boundary values $g\in C^{2+\alpha}(\partial\Omega,\mathbb R)$ there exists
a unique solution $u$ of the Dirichlet problem (\ref{dirichletproblem}).
\end{theorem}
{\it Proof:} \\
The uniqueness part follows from the assumption $\mathcal H_z\geq 0$ together 
with the maximum principle. For the existence part, 
consider a parameter $t\in [0,1]$ and the family of Dirichlet problems
\bee\label{l60}
u\in C^{2+\alpha}(\ol\Omega,\mathbb R) \quad , \quad
\sum_{i,j=1}^n G_{ij}(-\nabla u,1) \partial_{ij} u=t\,\mathcal H(x,u) \quad \mbox{in}\; \Omega
\quad \mbox{and} \quad  u=t\, g \quad \mbox{on}\; \partial\Omega \; . 
\ee
Because $G$ and $\mathcal H$ are assumed to be $C^{1+\alpha}$-functions, any such solution $u$
will belong to $C^{3+\alpha}(\Omega,\mathbb R)$ by interior Schauder theory.
Then the $C^0$-estimate Theorem \ref{t1}, the boundary gradient estimate Theorem \ref{t_boundarygradient}
together with the gradient maximum principle Corollary \ref{c1} 
yield a constant $C$ independent of $t$ such that
$$ ||u||_{C^1(\Omega)}\leq C $$
for any solution $u=u_t$ of (\ref{l60}). Using the Leray-Schauder Theorem \cite[Theorem 13.8]{gilbarg}
the Dirichlet problem (\ref{l60}) is solvable for any $t\in [0,1]$. For $t=1$ 
we obtain the desired solution of (\ref{dirichletproblem}). \hfill $\Box$ \\ \\
{\it Remarks:}
\begin{itemize}
\item[1)] The existence theorem applies in particular to the case $\mathcal H\equiv 0$
where we obtain graphs with vanishing weighted mean curvature, i.e. weighted minimal graphs.
Note that the assumptions (\ref{a1}) and (\ref{a3}) are satisfied in this case while
assumption (\ref{a2}) reduces to $H_G^-(x)<0<H_G^+(x)$ for $x\in\partial\Omega$.
\item[2)] The first two assumptions (\ref{a1}) and (\ref{a2}) in c) are natural in the sence that similar
assumptions are also needed for the classical prescribed mean curvature equation. The necessity the smallness
assumption (\ref{a1}), needed to obtain a $C^0$-estimate, is demonstrated by Theorem~\ref{nonexistence}.
The necessity of the boundary curvature condition (\ref{a2}), needed to obtain a boundary gradient estimate,
is demonstrated by the non-existence theorem \cite[Theorem 14.12]{gilbarg}
\item[3)] Assumption (\ref{a3}), required to obtain a maximum principle for the gradient,
may be relaxed somehow. Note however that in case of constant weighted mean curvature
$\mathcal H(x,z)\equiv h\in\mathbb R$ this assumption is satisfied.
Graphs of constant weighted mean curvature are of geometric interest as they arise
as critical points of the functional
$$ A(u):=\int_\Omega\Big( F(-\nabla u,1)+h u\Big) dx $$
(see Example \ref{example1}).
Considering $h\in\mathbb R$ as a Lagrange parameter, one looks for critical points
of $\int_\Omega F(-\nabla u,1)dx$ under the volume constraint $\int_\Omega u dx=\mbox{const}$.
\end{itemize}

Matthias Bergner, Jens Dittrich\\
Universit\"at Ulm, Fakult\"at f\"ur Mathematik und Wirtschaftswissenschaften, Institut f\"ur Analysis\\
Helmholtzstr. 18, D-89069 Ulm, Germany\\[0.3cm]
e-mail: matthias.bergner@uni-ulm.de, jens.dittrich@uni-ulm.de\\

\begin{thebibliography}{10}
\bibitem{blaschke}
W.Blaschke, K.Leichtwei\ss: {\it Elementare Differentialgeometrie}
Springer Berlin Heidelberg New York, 1973.
\bibitem{clarenz3}
U.Clarenz, H.von der Mosel:{\it Compactness theorems and an isoperimetric inequality
for critical points of elliptic parametric functionals}. Calc. Var. 12, 85--107, 2001.
\bibitem{clarenz1}
U.Clarenz: {\it Enclosure theorems for extremals of elliptic parametric functionals}.
Calc. Var. 15, 313--324, 2002.
\bibitem{clarenz2}
U.Clarenz, H.von der Mosel: {On surfaces of prescribed $F$-mean curvature}.
Pacific J. Math. 213, No. 1, 2004.
\bibitem{dubrovin}
B.A.Dubrovin, A.T.Fomenko, S.P.Novikov: {\it Modern Geometry-Methods and Applications}.
Graduate Texts in Mathematics, Springer Berlin Heidelberg New York, 1984.
\bibitem{froehlich}
S.Fr\"ohlich: {\it On two-dimensional immersions that are stable for parametric
functionals of constant mean curvature type}. 
Differential Geometry and its Applications 23, 235--256, 2005.
\bibitem{gilbarg}
D.Gilbarg, N.S.Trudinger: {\it Elliptic Partial Differential Equations of Second Order}. 
Springer, Berlin Heidelberg New York, 1983.
\bibitem{saudok}
F.Sauvigny: {\it Fl\"achen vorgeschriebener mittlerer Kr\"ummung mit
eineindeutiger Projektion auf eine Ebene}. 
Math. Zeit. 180, 41--67, 1982.
\bibitem{saubuch3}
F.Sauvigny: {\it Partial Differential Equations, Vol. 1 and 2}. Springer Universitext, 2006.
\end{thebibliography}
\end{document}